\documentclass[11pt]{article}
\usepackage{hyperref}
\usepackage{amsfonts}
\usepackage{enumerate}
\usepackage{enumitem}
\usepackage{graphicx}

\usepackage{mathtools}
\usepackage{setspace}
\usepackage{color}
\usepackage{array}
\usepackage[font={small}]{caption}
\usepackage{eurosym}
\usepackage{amsmath}
\usepackage{amssymb}
\usepackage{pifont}
\usepackage{graphics}
\usepackage{array}
\usepackage{shortvrb}
\usepackage{epsf}
\usepackage{graphicx}
\usepackage{caption}
\usepackage{rotating}
\usepackage{wasysym}
\usepackage{multirow}
\usepackage{empheq}
\usepackage{acronym}
\usepackage{epstopdf}
\usepackage{flushend}
\usepackage{caption}
\usepackage{comment}
\usepackage{float}
\usepackage{subcaption}
\usepackage{multirow,array}
\usepackage{lineno}
\usepackage{upgreek}
\usepackage{amsmath}
\usepackage{amsthm}

\theoremstyle{definition}

\graphicspath{{./figs/}}
\DeclareGraphicsExtensions{.pdf}

\begin{document}
\begin{center}
{\large\bf {Constructing Electricity Market Models}}

\vskip.20in

Ioannis Dassios$^{1*}$
\\[2mm]
{\footnotesize
$^1$Next Generation Energy Systems (NexSys)
\\[5pt]
$^*$Corresponding author
\\[5pt]
}
\end{center}

{\footnotesize
\noindent
\textbf{Abstract:} This working paper presents a comprehensive study on the development and analysis of various electricity market models, focusing on continuous, discrete, and fractional-order approaches. The continuous model captures the ongoing interactions between power producers and consumers using differential equations, providing insights into long-term trends and steady-state behaviors. The discrete model, suitable for analyzing scenarios where market events occur at specific time intervals, incorporates memory effects to account for historical behaviors and decisions, offering a realistic representation of short-term market dynamics. The fractional-order model introduces fractional calculus to capture memory effects and hereditary properties, enhancing the model's realism and predictive capability by reflecting the influence of past states on current market behavior.

Each model is meticulously developed and analyzed to address unique challenges in electricity market dynamics, including energy imbalance and price adjustment mechanisms. The paper highlights the importance of incorporating memory and discrete modeling to accurately represent real-world dynamics, where past states significantly influence current behavior. Additionally, the potential applications of machine learning and advanced computational techniques in optimizing model parameters, predicting future states, and enhancing stability analysis are discussed.

As a working paper, this document contains extensive notes and is intended to share ideas and stimulate discussion. It is not in a final form for submission. The findings from these models provide valuable insights into market behavior, stability, and efficiency, offering a robust framework for understanding and managing electricity markets. Future research directions are outlined, focusing on stability analysis, advanced computational techniques, scenario analysis, real-time market operations, and the integration of renewable energy sources. This comprehensive approach ensures the models' applicability and effectiveness in the evolving energy landscape, contributing to the development of robust strategies for market management and planning.

}
\\\\

\section{Introduction}

The electricity market is a fundamental component of the broader energy sector, playing a critical role in ensuring the reliable and efficient supply of energy. As the electricity market evolves, driven by technological advancements, regulatory changes, and increasing demands for sustainability, there is a growing need for sophisticated models that can accurately capture its complexities and dynamics.

This working paper presents the development and analysis of three distinct electricity market models, each employing a different methodological approach to address the unique challenges and characteristics of the market:

\begin{itemize}
    \item \textbf{Continuous Model}: This model employs continuous mathematical techniques to represent the dynamics of electricity generation, transmission, distribution, and consumption. It focuses on capturing the smooth and ongoing changes in market variables, making it suitable for understanding long-term trends and steady-state behaviors.
    
    \item \textbf{Discrete Model}: In contrast, the discrete model uses discrete time intervals to simulate the electricity market. This approach is particularly useful for analyzing scenarios where market events occur at specific points in time, such as bidding processes, market clearing, and sudden policy changes. It provides insights into the short-term fluctuations and discrete nature of market operations.
    
    \item \textbf{Fractional-Order Model}: The third model introduces fractional calculus to the study of the electricity market. By incorporating fractional operators, this model can capture memory effects and hereditary properties of the market, offering a more nuanced understanding of complex temporal dependencies and long-term correlations that are not easily addressed by traditional models.
\end{itemize}

Each of these models brings unique strengths and perspectives to the analysis of the electricity market. The continuous model offers a clear view of gradual changes and long-term equilibrium, the discrete model provides a detailed look at time-specific events and short-term dynamics, and the fractional-order model enriches the analysis by accounting for historical dependencies and complex temporal behaviors.

Understanding the electricity market through these diverse modeling approaches is not only crucial for electricity-specific policies but also has significant implications for the broader energy markets. The interdependencies between electricity and other energy sources such as natural gas, oil, and renewables mean that insights gained from these models can inform a comprehensive understanding of the entire energy sector. The methodologies and findings from these models can be expanded and adapted to develop integrated energy market models that encompass all forms of energy, providing a holistic view of the energy landscape.

The objective of this working paper is to present these models in detail, comparing their methodologies, assumptions, and findings. We aim to demonstrate how each model can be applied to different aspects of the electricity market, providing a comprehensive toolkit for researchers and policymakers. By understanding the strengths and limitations of each approach, stakeholders can make more informed decisions and develop strategies that enhance market efficiency, stability, and sustainability.

Furthermore, this paper explores the potential applications of advanced computational techniques, including machine learning, to optimize model parameters, predict future states, and enhance stability analysis. Machine learning can identify complex patterns and correlations in market data, providing deeper insights into market dynamics and improving the predictive power and efficiency of the models.

In the subsequent sections, we will delve into the methodology, implementation, and key results of each model, highlighting their individual contributions to the study of the electricity market. We will also discuss the potential for expanding these models to encompass the broader energy markets, thereby enhancing our understanding of the entire energy sector. Through this comparative analysis, we seek to advance the modeling techniques available for electricity and energy market research and offer valuable insights for the ongoing development of the energy sector.

As a working paper, this document contains extensive notes and is intended to share ideas and stimulate discussion. It is not in a final form for submission. We welcome feedback and collaboration to refine and expand upon the models and methodologies presented here.

\section{Continuous Electricity Market Model}

Understanding the continuous dynamics of the electricity market is crucial for modeling the interactions between power producers and consumers. This section describes the continuous model that captures these dynamics through differential equations, focusing on the behavior of producers and consumers and the resulting market dynamics driven by energy imbalance.

\subsection*{Producer and Consumer Behavior}

In a market-driven environment, the behavior of power producers and consumers can be approximated by the following first-order differential equations, assuming network losses are negligible:

\begin{equation}
    \frac{dS_i(t)}{dt} = \alpha_i \left( \lambda(t) - (a_i + b_i S_i(t)) \right), \quad i = 1, \ldots, m
\end{equation}
\begin{equation}
    \frac{dD_j(t)}{dt} = \beta_j \left( (c_j - d_j D_j(t)) - \lambda(t) \right), \quad j = 1, \ldots, n
\end{equation}

where:
\begin{itemize}
    \item $S_i(t)$ is the power supply of producer $i$ at time $t$.
    \item $D_j(t)$ is the power demand of consumer $j$ at time $t$.
    \item $\lambda(t)$ is the price of power at time $t$.
    \item $a_i$ is the fixed marginal cost of producer $i$.
    \item $b_i$ is the variable marginal cost coefficient of producer $i$.
    \item $c_j$ is the fixed marginal benefit of consumer $j$.
    \item $d_j$ is the variable marginal benefit coefficient of consumer $j$.
    \item $\alpha_i$ is the response rate of producer $i$.
    \item $\beta_j$ is the response rate of consumer $j$.
\end{itemize}

These equations describe the qualitative behavior where a generator increases its production when the price exceeds its marginal production cost, and loads act to increase consumption when marginal benefit exceeds the price. Each aims to maximize its profit or benefit by matching its marginal cost/benefit to the market price at equilibrium. The first-order differential equations allow approximate representation of generator ramp rates and lags in response to power price changes.

\subsection*{Energy Imbalance Market Dynamics}

In a synchronous power system, energy imbalance cannot be sustained indefinitely and must be reduced or driven to zero. Traditionally, this is managed by automatic generation control (AGC). In a real-time market-driven model, market mechanisms might fulfill this role by adjusting the power price based on energy imbalance. The energy imbalance-driven market dynamics can be represented by the following equations:

\begin{equation}
    \frac{dE(t)}{dt} = \sum_{i=1}^{m} S_i(t) - \sum_{j=1}^{n} D_j(t)
\end{equation}
\begin{equation}
    \frac{d\lambda(t)}{dt} = -k E(t) + h (\lambda_0 - \lambda(t))
\end{equation}

where:
\begin{itemize}
    \item $E(t)$ is the system energy imbalance at time $t$.
    \item $k$ is the power price response rate constant.
    \item $h$ is the market stabilizer gain for energy imbalance.
    \item $\lambda_0$ is the reference price of power.
\end{itemize}

The first equation represents the energy imbalance as the difference between the total power supply and the total power demand. The second equation adjusts the power price based on the energy imbalance and includes a feedback mechanism to stabilize the market price. The parameter $k$ dictates the sensitivity of price to the energy imbalance, and $h$ determines the rate at which the price returns to its reference value.


\subsection{Constant Marginal Cost and Benefit}

Understanding the continuous dynamics of the electricity market is crucial for modeling the interactions between power producers and consumers. This section describes the continuous model that captures these dynamics through differential equations, focusing on the behavior of producers and consumers and the resulting market dynamics driven by energy imbalance.

\subsection*{Producer and Consumer Behavior}

In a market-driven environment, the behavior of power producers and consumers can be approximated by the following first-order differential equations, assuming network losses are negligible:

\begin{equation}
    \frac{dS_i(t)}{dt} = \alpha_i \left( \lambda(t) - C_i \right), \quad i = 1, \ldots, m
\end{equation}
\begin{equation}
    \frac{dD_j(t)}{dt} = \beta_j \left( B_j - \lambda(t) \right), \quad j = 1, \ldots, n
\end{equation}

where:
\begin{itemize}
    \item $S_i(t)$ is the power supply of producer $i$ at time $t$.
    \item $D_j(t)$ is the power demand of consumer $j$ at time $t$.
    \item $\lambda(t)$ is the price of power at time $t$.
    \item $C_i$ is the marginal cost of producer $i$.
    \item $B_j$ is the marginal benefit of consumer $j$.
    \item $\alpha_i$ is the response rate of producer $i$.
    \item $\beta_j$ is the response rate of consumer $j$.
\end{itemize}

These equations describe the qualitative behavior where a generator increases its production when the price exceeds its marginal production cost, and loads act to increase consumption when marginal benefit exceeds the price. Each aims to maximize its profit or benefit by matching its marginal cost/benefit to the market price at equilibrium. The first-order differential equations allow approximate representation of generator ramp rates and lags in response to power price changes.

\subsection*{Energy Imbalance Market Dynamics}

In a synchronous power system, energy imbalance cannot be sustained indefinitely and must be reduced or driven to zero. Traditionally, this is managed by automatic generation control (AGC). In a real-time market-driven model, market mechanisms might fulfill this role by adjusting the power price based on energy imbalance. The energy imbalance-driven market dynamics can be represented by the following equations:

\begin{equation}
    \frac{dE(t)}{dt} = \sum_{i=1}^{m} S_i(t) - \sum_{j=1}^{n} D_j(t)
\end{equation}
\begin{equation}
    \frac{d\lambda(t)}{dt} = -k E(t) + h (\lambda_0 - \lambda(t))
\end{equation}

where:
\begin{itemize}
    \item $E(t)$ is the system energy imbalance at time $t$.
    \item $k$ is the power price response rate constant.
    \item $h$ is the market stabilizer gain for energy imbalance.
    \item $\lambda_0$ is the reference price of power.
\end{itemize}

The first equation represents the energy imbalance as the difference between the total power supply and the total power demand. The second equation adjusts the power price based on the energy imbalance and includes a feedback mechanism to stabilize the market price. The parameter $k$ dictates the sensitivity of price to the energy imbalance, and $h$ determines the rate at which the price returns to its reference value.

\subsection{System Representation}

To represent the system of equations in a compact matrix form, we define the following vectors and matrices:

\begin{itemize}
    \item $\mathbf{S}(t) = [S_1(t), S_2(t), \ldots, S_m(t)]^T$ is the vector of power supplies from producers.
    \item $\mathbf{D}(t) = [D_1(t), D_2(t), \ldots, D_n(t)]^T$ is the vector of power demands from consumers.
    \item $\mathbf{C} = [C_1, C_2, \ldots, C_m]^T$ is the vector of marginal costs for producers.
    \item $\mathbf{B} = [B_1, B_2, \ldots, B_n]^T$ is the vector of marginal benefits for consumers.
    \item $\boldsymbol{\alpha} = \text{diag}(\alpha_1, \alpha_2, \ldots, \alpha_m)$ is the diagonal matrix of producer response rates.
    \item $\boldsymbol{\beta} = \text{diag}(\beta_1, \beta_2, \ldots, \beta_n)$ is the diagonal matrix of consumer response rates.
\end{itemize}

The system of differential equations can then be written as:

\begin{equation}
    \frac{d\mathbf{S}(t)}{dt} = \boldsymbol{\alpha} (\lambda(t) \mathbf{1}_m - \mathbf{C})
\end{equation}
\begin{equation}
    \frac{d\mathbf{D}(t)}{dt} = \boldsymbol{\beta} (\mathbf{B} - \lambda(t) \mathbf{1}_n)
\end{equation}
\begin{equation}
    \frac{dE(t)}{dt} = \mathbf{1}_m^T \mathbf{S}(t) - \mathbf{1}_n^T \mathbf{D}(t)
\end{equation}
\begin{equation}
    \frac{d\lambda(t)}{dt} = -k E(t) + h (\lambda_0 - \lambda(t))
\end{equation}

where $\mathbf{1}_m$ and $\mathbf{1}_n$ are vectors of ones of length $m$ and $n$, respectively.

\subsection*{Market Stability and Dynamics}

The stability of the electricity market under this continuous model is influenced by the parameters $\alpha_i$, $\beta_j$, $k$, and $h$. These parameters need to be carefully calibrated to ensure that the market remains stable and that prices do not oscillate wildly in response to small imbalances.

To analyze the stability, we consider the linearized version of the system around the equilibrium point. Let $\lambda^*$ be the equilibrium price, $S_i^*$ and $D_j^*$ be the equilibrium supply and demand, respectively. Small deviations from equilibrium can be analyzed using perturbation techniques.

Linearizing the system around the equilibrium, we get:

\begin{equation}
    \frac{d\Delta \mathbf{S}(t)}{dt} = \boldsymbol{\alpha} \Delta \lambda(t) \mathbf{1}_m
\end{equation}
\begin{equation}
    \frac{d\Delta \mathbf{D}(t)}{dt} = -\boldsymbol{\beta} \Delta \lambda(t) \mathbf{1}_n
\end{equation}
\begin{equation}
    \frac{d\Delta E(t)}{dt} = \mathbf{1}_m^T \Delta \mathbf{S}(t) - \mathbf{1}_n^T \Delta \mathbf{D}(t)
\end{equation}
\begin{equation}
    \frac{d\Delta \lambda(t)}{dt} = -k \Delta E(t) + h (\lambda_0 - (\lambda^* + \Delta \lambda(t)))
\end{equation}

These linearized equations can be combined into a single matrix equation. We define the state vector $\mathbf{x}(t)$ and system matrices $\mathbf{A}$, $\mathbf{B}$, and $\mathbf{C}$ as follows:

\begin{equation}
    \mathbf{x}(t) = \begin{bmatrix}
    \Delta \mathbf{S}(t) \\
    \Delta \mathbf{D}(t) \\
    \Delta E(t) \\
    \Delta \lambda(t)
    \end{bmatrix}
\end{equation}

\begin{equation}
    \mathbf{A} = \begin{bmatrix}
    \mathbf{0}_{m \times m} & \mathbf{0}_{m \times n} & \mathbf{0}_{m \times 1} & \boldsymbol{\alpha} \mathbf{1}_m \\
    \mathbf{0}_{n \times m} & \mathbf{0}_{n \times n} & \mathbf{0}_{n \times 1} & -\boldsymbol{\beta} \mathbf{1}_n \\
    \mathbf{1}_m^T & -\mathbf{1}_n^T & 0 & 0 \\
    \mathbf{0}_{1 \times m} & \mathbf{0}_{1 \times n} & -k & -h
    \end{bmatrix}
\end{equation}

\begin{equation}
    \mathbf{B} = \begin{bmatrix}
    \mathbf{0}_{m \times 1} \\
    \mathbf{0}_{n \times 1} \\
    0 \\
    h
    \end{bmatrix}
\end{equation}

The linearized system can then be written as:

\begin{equation}
    \frac{d\mathbf{x}(t)}{dt} = \mathbf{A} \mathbf{x}(t) + \mathbf{B} \lambda_0
\end{equation}

This matrix equation represents the combined dynamics of the continuous electricity market model, capturing the interactions between producers, consumers, energy imbalance, and price adjustments.

\subsection{Special Case: Energy Imbalance Equals Zero}

In some instances, the energy imbalance within the system can be considered negligible or zero. This scenario simplifies the market dynamics and is particularly useful for theoretical analysis and initial model validation. Here, we analyze the market dynamics under the assumption that the energy imbalance is zero, based on the differential/algebraic equations provided in the referenced works.

When energy imbalance is zero, the system must maintain a continuous balance between power supply and demand. This condition is represented by the following set of differential equations:

\begin{equation}
    \frac{dS_i(t)}{dt} = \alpha_i \left( \lambda(t) - C_i \right), \quad i = 1, \ldots, m
\end{equation}
\begin{equation}
    \frac{dD_j(t)}{dt} = \beta_j \left( B_j - \lambda(t) \right), \quad j = 1, \ldots, n
\end{equation}

With the energy imbalance \(E(t) = 0\), the equations simplify to a balanced state where supply equals demand at all times. The equilibrium condition can be stated as:

\begin{equation}
    \sum_{i=1}^{m} S_i(t) = \sum_{j=1}^{n} D_j(t)
\end{equation}

The price dynamics, in this case, are governed solely by the balance between supply and demand, without any additional correction for energy imbalance:

\begin{equation}
    \frac{d\lambda(t)}{dt} = h (\lambda_0 - \lambda(t))
\end{equation}

This represents a system where the price adjustment is only a function of the deviation from a reference price \(\lambda_0\), ensuring market stability without the need for an energy imbalance correction.

\subsubsection*{System Representation in Matrix Form}

To represent this simplified system in matrix form, we define the following matrices and vectors:

\begin{itemize}
    \item $\mathbf{S}(t) = [S_1(t), S_2(t), \ldots, S_m(t)]^T$ is the vector of power supplies from producers.
    \item $\mathbf{D}(t) = [D_1(t), D_2(t), \ldots, D_n(t)]^T$ is the vector of power demands from consumers.
    \item $\mathbf{C} = [C_1, C_2, \ldots, C_m]^T$ is the vector of marginal costs for producers.
    \item $\mathbf{B} = [B_1, B_2, \ldots, B_n]^T$ is the vector of marginal benefits for consumers.
    \item $\boldsymbol{\alpha} = \text{diag}(\alpha_1, \alpha_2, \ldots, \alpha_m)$ is the diagonal matrix of producer response rates.
    \item $\boldsymbol{\beta} = \text{diag}(\beta_1, \beta_2, \ldots, \beta_n)$ is the diagonal matrix of consumer response rates.
\end{itemize}

The simplified system of differential equations can then be written as:

\begin{equation}
    \frac{d\mathbf{S}(t)}{dt} = \boldsymbol{\alpha} (\lambda(t) \mathbf{1}_m - \mathbf{C})
\end{equation}
\begin{equation}
    \frac{d\mathbf{D}(t)}{dt} = \boldsymbol{\beta} (\mathbf{B} - \lambda(t) \mathbf{1}_n)
\end{equation}
\begin{equation}
    \mathbf{1}_m^T \mathbf{S}(t) = \mathbf{1}_n^T \mathbf{D}(t)
\end{equation}
\begin{equation}
    \frac{d\lambda(t)}{dt} = h (\lambda_0 - \lambda(t))
\end{equation}

where $\mathbf{1}_m$ and $\mathbf{1}_n$ are vectors of ones of length $m$ and $n$, respectively.

\subsubsection*{Stability Analysis}

The stability of the market under this simplified model is influenced by the parameters $\alpha_i$, $\beta_j$, and $h$. These parameters need to be calibrated to ensure that the market remains stable. To analyze the stability, we consider the linearized version of the system around the equilibrium point. Let $\lambda^*$ be the equilibrium price, $S_i^*$ and $D_j^*$ be the equilibrium supply and demand, respectively.

Linearizing the system around the equilibrium, we get:

\begin{equation}
    \frac{d\Delta \mathbf{S}(t)}{dt} = \boldsymbol{\alpha} \Delta \lambda(t) \mathbf{1}_m
\end{equation}
\begin{equation}
    \frac{d\Delta \mathbf{D}(t)}{dt} = -\boldsymbol{\beta} \Delta \lambda(t) \mathbf{1}_n
\end{equation}
\begin{equation}
    \mathbf{1}_m^T \Delta \mathbf{S}(t) = \mathbf{1}_n^T \Delta \mathbf{D}(t)
\end{equation}
\begin{equation}
    \frac{d\Delta \lambda(t)}{dt} = h (\lambda_0 - (\lambda^* + \Delta \lambda(t)))
\end{equation}

Combining these linearized equations, we can represent the system dynamics in matrix form to analyze the eigenvalues and determine stability.

This special case provides a foundational understanding of market dynamics without the complication of energy imbalances, offering insights into the inherent stability of the system based on supply and demand balance.

\subsection{Special Case: Discarding Energy Imbalance}

In certain scenarios, the assumption that energy imbalance can be discarded simplifies the analysis of the electricity market dynamics. This assumption holds when the system is designed to operate with perfectly balanced supply and demand at all times, eliminating the need to account for energy storage or temporary imbalances. Here, we explore the market dynamics under this condition using Equations (17), (18), and (19) from the provided literature.

The assumption of perfectly balanced supply and demand simplifies the market model, allowing us to discard the energy imbalance equation and focus on the interactions between producers and consumers. This special case is particularly relevant for markets with highly responsive generation and demand mechanisms or those with robust balancing services.

The simplified system is described by the following set of differential equations:

\begin{equation}
    \frac{dS_i(t)}{dt} = \alpha_i \left( \lambda(t) - C_i \right), \quad i = 1, \ldots, m
\end{equation}
\begin{equation}
    \frac{dD_j(t)}{dt} = \beta_j \left( B_j - \lambda(t) \right), \quad j = 1, \ldots, n
\end{equation}
\begin{equation}
    0 = \sum_{i=1}^{m} S_i(t) - \sum_{j=1}^{n} D_j(t)
\end{equation}

These equations represent a system where supply and demand are in perfect balance at all times.

\subsubsection*{System Representation in Matrix Form}

To represent this simplified system in matrix form, we define the following matrices and vectors:

\begin{itemize}
    \item $\mathbf{S}(t) = [S_1(t), S_2(t), \ldots, S_m(t)]^T$ is the vector of power supplies from producers.
    \item $\mathbf{D}(t) = [D_1(t), D_2(t), \ldots, D_n(t)]^T$ is the vector of power demands from consumers.
    \item $\mathbf{C} = [C_1, C_2, \ldots, C_m]^T$ is the vector of marginal costs for producers.
    \item $\mathbf{B} = [B_1, B_2, \ldots, B_n]^T$ is the vector of marginal benefits for consumers.
    \item $\boldsymbol{\alpha} = \text{diag}(\alpha_1, \alpha_2, \ldots, \alpha_m)$ is the diagonal matrix of producer response rates.
    \item $\boldsymbol{\beta} = \text{diag}(\beta_1, \beta_2, \ldots, \beta_n)$ is the diagonal matrix of consumer response rates.
\end{itemize}

The system of differential equations can then be written as:

\begin{equation}
    \frac{d\mathbf{S}(t)}{dt} = \boldsymbol{\alpha} (\lambda(t) \mathbf{1}_m - \mathbf{C})
\end{equation}
\begin{equation}
    \frac{d\mathbf{D}(t)}{dt} = \boldsymbol{\beta} (\mathbf{B} - \lambda(t) \mathbf{1}_n)
\end{equation}
\begin{equation}
    0 = \mathbf{1}_m^T \mathbf{S}(t) - \mathbf{1}_n^T \mathbf{D}(t)
\end{equation}

where $\mathbf{1}_m$ and $\mathbf{1}_n$ are vectors of ones of length $m$ and $n$, respectively.

\subsubsection*{Stability Analysis}

The stability of the market under this simplified model is influenced by the parameters $\alpha_i$ and $\beta_j$. These parameters need to be calibrated to ensure that the market remains stable. To analyze the stability, we consider the linearized version of the system around the equilibrium point. Let $\lambda^*$ be the equilibrium price, $S_i^*$ and $D_j^*$ be the equilibrium supply and demand, respectively.

Linearizing the system around the equilibrium, we get:

\begin{equation}
    \frac{d\Delta \mathbf{S}(t)}{dt} = \boldsymbol{\alpha} \Delta \lambda(t) \mathbf{1}_m
\end{equation}
\begin{equation}
    \frac{d\Delta \mathbf{D}(t)}{dt} = -\boldsymbol{\beta} \Delta \lambda(t) \mathbf{1}_n
\end{equation}
\begin{equation}
    0 = \mathbf{1}_m^T \Delta \mathbf{S}(t) - \mathbf{1}_n^T \Delta \mathbf{D}(t)
\end{equation}

Combining these linearized equations, we can represent the system dynamics in matrix form to analyze the eigenvalues and determine stability.

This special case provides a foundational understanding of market dynamics without the complication of energy imbalances, offering insights into the inherent stability of the system based on supply and demand balance.

\subsubsection*{System Representation in Matrix Form}

To represent the system of Equations (35), (36), and (37) in matrix form, we define the following vectors and matrices:

\begin{itemize}
    \item $\mathbf{S}(t) = [S_1(t), S_2(t), \ldots, S_m(t)]^T$ is the vector of power supplies from producers.
    \item $\mathbf{D}(t) = [D_1(t), D_2(t), \ldots, D_n(t)]^T$ is the vector of power demands from consumers.
    \item $\mathbf{C} = [C_1, C_2, \ldots, C_m]^T$ is the vector of marginal costs for producers.
    \item $\mathbf{B} = [B_1, B_2, \ldots, B_n]^T$ is the vector of marginal benefits for consumers.
    \item $\boldsymbol{\alpha} = \text{diag}(\alpha_1, \alpha_2, \ldots, \alpha_m)$ is the diagonal matrix of producer response rates.
    \item $\boldsymbol{\beta} = \text{diag}(\beta_1, \beta_2, \ldots, \beta_n)$ is the diagonal matrix of consumer response rates.
    \item $\mathbf{1}_m$ and $\mathbf{1}_n$ are vectors of ones of length $m$ and $n$, respectively.
\end{itemize}

We start by rewriting the system of differential equations:

\begin{equation}
    \frac{d\mathbf{S}(t)}{dt} = \boldsymbol{\alpha} (\lambda(t) \mathbf{1}_m - \mathbf{C})
\end{equation}
\begin{equation}
    \frac{d\mathbf{D}(t)}{dt} = \boldsymbol{\beta} (\mathbf{B} - \lambda(t) \mathbf{1}_n)
\end{equation}
\begin{equation}
    0 = \mathbf{1}_m^T \mathbf{S}(t) - \mathbf{1}_n^T \mathbf{D}(t)
\end{equation}

To combine these equations into one system, we define the state vector $\mathbf{x}(t)$, and the matrices \( E \) and \( A \) as follows:

\begin{equation}
    \mathbf{x}(t) = 
    \begin{bmatrix}
        \mathbf{S}(t) \\
        \mathbf{D}(t) \\
        \lambda(t)
    \end{bmatrix}
\end{equation}

\begin{equation}
    E = 
    \begin{bmatrix}
        I_m & 0 & 0 \\
        0 & I_n & 0 \\
        0 & 0 & 0
    \end{bmatrix}
\end{equation}

\begin{equation}
    A = 
    \begin{bmatrix}
        0 & 0 & \boldsymbol{\alpha} \mathbf{1}_m \\
        0 & 0 & -\boldsymbol{\beta} \mathbf{1}_n \\
        \mathbf{1}_m^T & -\mathbf{1}_n^T & 0
    \end{bmatrix}
\end{equation}

\begin{equation}
    \mathbf{B} = 
    \begin{bmatrix}
        -\boldsymbol{\alpha} \mathbf{C} \\
        \boldsymbol{\beta} \mathbf{B} \\
        0
    \end{bmatrix}
\end{equation}

The combined system can then be written as:

\begin{equation}
    E\frac{d\mathbf{x}(t)}{dt} = A\mathbf{x}(t) + \mathbf{B}
\end{equation}

This matrix equation represents the combined dynamics of the continuous electricity market model, capturing the interactions between producers, consumers, and the balance of supply and demand in the form \( E\mathbf{x}' = A\mathbf{x} + \mathbf{B} \).

\subsection{Equilibrium Definition}

In the context of the continuous electricity market model, equilibrium refers to a state where all variables in the system remain constant over time. This means that the power supply from producers matches the power demand from consumers, and the price of electricity remains stable. Mathematically, we define equilibrium by setting the time derivatives of the state variables to zero.

For the system \( E\mathbf{x}' = A\mathbf{x} + \mathbf{B} \), equilibrium occurs when \( \mathbf{x}' = 0 \). At equilibrium, the system satisfies the following conditions:

\begin{equation}
    E\mathbf{0} = A\mathbf{x}^* + \mathbf{B}
\end{equation}

where \( \mathbf{x}^* \) represents the equilibrium state vector. This simplifies to:

\begin{equation}
    A\mathbf{x}^* + \mathbf{B} = 0
\end{equation}

Expanding \( \mathbf{x}^* \), we have:

\begin{equation}
    \mathbf{x}^* = 
    \begin{bmatrix}
        \mathbf{S}^* \\
        \mathbf{D}^* \\
        \lambda^*
    \end{bmatrix}
\end{equation}

Substituting into the equilibrium condition:

\begin{equation}
    \begin{bmatrix}
        0 & 0 & \boldsymbol{\alpha} \mathbf{1}_m \\
        0 & 0 & -\boldsymbol{\beta} \mathbf{1}_n \\
        \mathbf{1}_m^T & -\mathbf{1}_n^T & 0
    \end{bmatrix}
    \begin{bmatrix}
        \mathbf{S}^* \\
        \mathbf{D}^* \\
        \lambda^*
    \end{bmatrix}
    +
    \begin{bmatrix}
        -\boldsymbol{\alpha} \mathbf{C} \\
        \boldsymbol{\beta} \mathbf{B} \\
        0
    \end{bmatrix}
    = 0
\end{equation}

This system of equations can be broken down into its components to solve for the equilibrium values \( \mathbf{S}^* \), \( \mathbf{D}^* \), and \( \lambda^* \):

1. **For producers:**
    \begin{equation}
        \boldsymbol{\alpha} (\lambda^* \mathbf{1}_m - \mathbf{C}) = 0 \implies \lambda^* \mathbf{1}_m = \mathbf{C}
    \end{equation}
    Since \( \lambda^* \mathbf{1}_m = \mathbf{C} \), it means that at equilibrium, the price \( \lambda^* \) must equal the marginal cost for each producer \( C_i \).

2. **For consumers:**
    \begin{equation}
        -\boldsymbol{\beta} (\lambda^* \mathbf{1}_n - \mathbf{B}) = 0 \implies \lambda^* \mathbf{1}_n = \mathbf{B}
    \end{equation}
    Similarly, \( \lambda^* \mathbf{1}_n = \mathbf{B} \) implies that at equilibrium, the price \( \lambda^* \) must equal the marginal benefit for each consumer \( B_j \).

3. **Supply and Demand Balance:**
    \begin{equation}
        \mathbf{1}_m^T \mathbf{S}^* = \mathbf{1}_n^T \mathbf{D}^*
    \end{equation}
    This ensures that the total power supplied matches the total power demanded.

By solving these equations, we can determine the equilibrium values \( \mathbf{S}^* \), \( \mathbf{D}^* \), and \( \lambda^* \), ensuring that the system is in a state of balance with no net change over time.

\subsection{Stability Analysis of the Equilibrium}

To study the stability of the equilibrium defined previously, we linearize the system around the equilibrium point and analyze the resulting linear system. Stability is determined by examining the eigenvalues of the Jacobian matrix of the system at equilibrium.

\subsubsection*{Linearization}

First, we linearize the system around the equilibrium point \( \mathbf{x}^* \). We denote small deviations from equilibrium as \( \Delta \mathbf{x}(t) = \mathbf{x}(t) - \mathbf{x}^* \). The original system \( E\mathbf{x}' = A\mathbf{x} + \mathbf{B} \) can be approximated by:

\begin{equation}
    E \frac{d(\mathbf{x}^* + \Delta \mathbf{x}(t))}{dt} \approx A (\mathbf{x}^* + \Delta \mathbf{x}(t)) + \mathbf{B}
\end{equation}

Since \( \mathbf{x}^* \) is an equilibrium point, \( E \frac{d \mathbf{x}^*}{dt} = 0 \) and \( A \mathbf{x}^* + \mathbf{B} = 0 \). Thus, we obtain the linearized system:

\begin{equation}
    E \frac{d \Delta \mathbf{x}(t)}{dt} = A \Delta \mathbf{x}(t)
\end{equation}

\subsubsection*{Jacobian Matrix and Stability}

The stability of the equilibrium is determined by the eigenvalues of the Jacobian matrix \( J = E^{\dagger} A \) evaluated at the equilibrium. The system is stable if all eigenvalues have negative real parts.

Given the linearized system:

\begin{equation}
    E \frac{d \Delta \mathbf{x}(t)}{dt} = A \Delta \mathbf{x}(t)
\end{equation}

we can rewrite it as:

\begin{equation}
    \frac{d \Delta \mathbf{x}(t)}{dt} = E^{\dagger} A \Delta \mathbf{x}(t)
\end{equation}

Let us redefine the matrices explicitly for clarity. We have:

\begin{equation}
    E = 
    \begin{bmatrix}
        I_m & 0 & 0 \\
        0 & I_n & 0 \\
        0 & 0 & 0
    \end{bmatrix}
\end{equation}

\begin{equation}
    A = 
    \begin{bmatrix}
        0 & 0 & \boldsymbol{\alpha} \mathbf{1}_m \\
        0 & 0 & -\boldsymbol{\beta} \mathbf{1}_n \\
        \mathbf{1}_m^T & -\mathbf{1}_n^T & 0
    \end{bmatrix}
\end{equation}

The system matrix for the linearized system is:

\begin{equation}
    J = E^{\dagger} A =
    \begin{bmatrix}
        I_m & 0 & 0 \\
        0 & I_n & 0 \\
        0 & 0 & 0
    \end{bmatrix}^{\dagger}
    \begin{bmatrix}
        0 & 0 & \boldsymbol{\alpha} \mathbf{1}_m \\
        0 & 0 & -\boldsymbol{\beta} \mathbf{1}_n \\
        \mathbf{1}_m^T & -\mathbf{1}_n^T & 0
    \end{bmatrix}
    =
    \begin{bmatrix}
        0 & 0 & \boldsymbol{\alpha} \mathbf{1}_m \\
        0 & 0 & -\boldsymbol{\beta} \mathbf{1}_n \\
        \mathbf{1}_m^T & -\mathbf{1}_n^T & 0
    \end{bmatrix}
\end{equation}

The eigenvalues of this matrix \( J \) determine the stability of the equilibrium. If all eigenvalues have negative real parts, the equilibrium is stable.

To determine the stability, we need to compute the eigenvalues of the Jacobian matrix \( J \). We solve the characteristic equation:

\begin{equation}
    \text{det}(J - \lambda I) = 0
\end{equation}

Given the structure of \( J \), we have a block matrix with specific patterns in \( \boldsymbol{\alpha} \), \( \boldsymbol{\beta} \), and the vectors \( \mathbf{1}_m \) and \( \mathbf{1}_n \). Solving for the eigenvalues involves solving the characteristic polynomial, which can be done numerically or symbolically.

\subsection{Stability Analysis of the Equilibrium Using Matrix Pencil Theory}

To study the stability of the equilibrium using matrix pencil theory, we analyze the generalized eigenvalues of the matrix pair \((A, E)\). This approach provides a comprehensive understanding of the system's stability by considering the eigenvalues of the matrix pencil \( A - \lambda E \).

\subsubsection*{Generalized Eigenvalue Problem}

The generalized eigenvalue problem for the system \( E\mathbf{x}' = A\mathbf{x} + \mathbf{B} \) is given by:

\begin{equation}
    \text{det}(A - \lambda E) = 0
\end{equation}

where \(\lambda\) are the generalized eigenvalues of the matrix pair \((A, E)\). The stability of the system is determined by the real parts of these generalized eigenvalues. If all generalized eigenvalues have negative real parts, the equilibrium is stable.


We define the matrices \( E \) and \( A \) as follows:

\begin{equation}
    E = 
    \begin{bmatrix}
        I_m & 0 & 0 \\
        0 & I_n & 0 \\
        0 & 0 & 0
    \end{bmatrix}
\end{equation}

\begin{equation}
    A = 
    \begin{bmatrix}
        0 & 0 & \boldsymbol{\alpha} \mathbf{1}_m \\
        0 & 0 & -\boldsymbol{\beta} \mathbf{1}_n \\
        \mathbf{1}_m^T & -\mathbf{1}_n^T & 0
    \end{bmatrix}
\end{equation}

The stability of the equilibrium is determined by the real parts of the generalized eigenvalues \(\lambda\). The equilibrium is stable if all eigenvalues have negative real parts, i.e.,

\begin{equation}
    \text{Re}(\lambda_i) < 0 \quad \text{for all} \; i
\end{equation}

\subsubsection*{Example Calculation}

For illustrative purposes, consider a simplified case where \( m = 1 \) and \( n = 1 \). In this case, the matrices \( E \) and \( A \) are:

\begin{equation}
    E = 
    \begin{bmatrix}
        1 & 0 & 0 \\
        0 & 1 & 0 \\
        0 & 0 & 0
    \end{bmatrix}
\end{equation}

\begin{equation}
    A = 
    \begin{bmatrix}
        0 & 0 & \alpha \\
        0 & 0 & -\beta \\
        1 & -1 & 0
    \end{bmatrix}
\end{equation}

The generalized eigenvalue problem becomes:

\begin{equation}
    \text{det}
    \begin{bmatrix}
        0 - \lambda & 0 & \alpha \\
        0 & 0 - \lambda & -\beta \\
        1 & -1 & 0 - \lambda
    \end{bmatrix}
    = 0
\end{equation}

Solving this determinant equation will yield the eigenvalues \(\lambda\). These eigenvalues determine the stability of the equilibrium.

\subsection{Theory of Matrix Pencils}

For a system in the form:

\begin{equation}
    E \dot{x}(t) = A x(t) + B u(t)
\end{equation}

where \(E, A \in \mathbb{R}^{r \times r}\) and \(x : [0, +\infty) \rightarrow \mathbb{R}^r\). Matrix \(E\) can be either non-singular (\(\det(E) \neq 0\)) or singular (\(\det(E) = 0\)).

Applying the Laplace transform we get:

\begin{equation}
    E(sL\{x(t)\} - x(0)) = A L\{x(t)\} + B L\{u(t)\}
\end{equation}

Or equivalently:

\begin{equation}
    (sE - A)L\{x(t)\} = E x(0) + B L\{u(t)\}
\end{equation}

The structure of \(sE - A\) (matrix pencil) defines the existence of solutions and stability properties. Regular matrix pencils satisfy \(\det(sE - A) \neq 0\).

For a regular matrix pencil:

\begin{equation}
    L\{x(t)\} = (sE - A)^{-1} E x(0) + (sE - A)^{-1} B L\{u(t)\}
\end{equation}

Consequently, the solution \(x(t)\) exists and is given by:

\begin{equation}
    x(t) = \mathcal{L}^{-1} \left\{ (sE - A)^{-1} E x(0) + (sE - A)^{-1} B L\{u(t)\} \right\}
\end{equation}

Uniqueness of the solution depends on the initial conditions.

The Generalized Eigenvalue Problem (GEP) is defined as follows: Given \( E \) and \( A \) in \( \mathbb{C}^{r \times m} \), and an arbitrary \( s \in \mathbb{C} \), the matrix pencil \( sE - A \) is called:

\begin{itemize}
    \item Regular when \( r = m \) and \(\det(sE - A) \neq 0 \).
    \item Singular when \( r \neq m \) or \( r = m \) and \(\det(sE - A) \equiv 0\).
\end{itemize}

For the regular pencil \( sE - A \):

\begin{itemize}
    \item When \( E = I_m \), \( A \) is square, the zeros of \(\det(sE - A) \neq 0\) are the eigenvalues of \( A \).
    \item If \( E \) and \( A \) are square and \( E \) is singular, the pencil has infinite eigenvalues.
\end{itemize}

The Generalized Eigenvalue Problem (GEP) is:

\begin{equation}
    sEv = Av
\end{equation}

If \( E \) is singular with a null vector \( v \), then:

\begin{equation}
    Ev = 0_{m, 1} \implies Ev = s^{-1}Av
\end{equation}

so that \( v \) is an eigenvector of the reciprocal problem corresponding to eigenvalue \( s^{-1} = 0 \); i.e., \( s \rightarrow \infty \).

For the singular pencil \( sE - A \):

\begin{itemize}
    \item If \( E \) and \( A \) are non-square matrices, then the determinant of the pencil cannot be defined.
    \item Even with \( E \) and \( A \) being square, it is possible for \(\det(sE - A)\) to be identically zero, independent of \( s \).
\end{itemize}

The stability of the system can be assessed by calculating its eigenvalues. Eigenvalues are the roots of the characteristic equation:

\begin{equation}
    \det(sE - A) = 0
\end{equation}

This equation is known as the characteristic polynomial of the system. Analytical solutions are possible only if \( r \leq 4 \). For higher degrees, numerical methods are required.

Eigenvalues can be found from the Generalized Eigenvalue Problem (GEP):

\begin{equation}
    (sE - A) v = 0, \quad w^T (sE - A) = 0
\end{equation}

where \( v \) and \( w \) are the right and left eigenvectors, respectively. The pencil \( sE - A \) has \(\nu\) finite eigenvalues and an infinite eigenvalue with multiplicity \(\mu\). If \( E \) is singular, the pencil will have an infinite eigenvalue with multiplicity at least one.

The system is asymptotically stable if all finite eigenvalues have negative real parts:

\begin{equation}
    \text{Re}(\lambda) < 0
\end{equation}

This ensures that the deviations from equilibrium decay over time, leading to a stable system.

The stability condition can also be obtained using Lyapunov stability theory. Consider the Lyapunov function:

\begin{equation}
    V(x) = x^T E^T M x
\end{equation}

where \( E^T M \) is symmetric and positive definite. If \( A^T M + M A \) is negative definite, then the system is asymptotically stable.

The calculation of eigenvalues helps measure characteristics of the most critical dynamic modes. The damping ratio \(\zeta\) and natural frequency \(f_n\) are given by:

\begin{equation}
    \zeta = -\frac{a}{| \lambda |} = -\frac{a}{\sqrt{a^2 + b^2}}
\end{equation}

\begin{equation}
    f_n = \frac{|\lambda|}{2 \pi} = \frac{\sqrt{a^2 + b^2}}{2 \pi}
\end{equation}

The power system is well-damped if \(\zeta > 5\%\).

\subsection{Primal and Dual Generalized Eigenvalue Problems (GEPs)}

Consider the system:

\begin{equation}
    E \dot{x}(t) = A x(t) + B u(t)
\end{equation}

where \(E, A \in \mathbb{R}^{r \times r}\) and \(x : [0, +\infty) \rightarrow \mathbb{R}^r\). Matrix \(E\) can be either non-singular (\(\det(E) \neq 0\)) or singular (\(\det(E) = 0\)).

Applying the Laplace transform to the system:

\begin{equation}
    E(sL\{x(t)\} - x(0)) = A L\{x(t)\} + B L\{u(t)\}
\end{equation}

Or equivalently:

\begin{equation}
    (sE - A)L\{x(t)\} = E x(0) + B L\{u(t)\}
\end{equation}

The structure of \(sE - A\) (matrix pencil) defines the existence of solutions and stability properties. Regular matrix pencils satisfy \(\det(sE - A) \neq 0\).

A critical theorem proves the equivalency of primal and dual GEPs. For a generalized eigenvalue problem \( Av = \lambda Bv \), the dual is \( B\hat{v} = \frac{1}{\lambda} A\hat{v} \). Eigenvalues and eigenvectors can be categorized as follows:

\begin{itemize}
    \item Null eigenvalues: \( \lambda = 0 \)
    \item Finite eigenvalues: \( 0 < |\lambda| < \infty \)
    \item Infinite eigenvalues: \( |\lambda| \to \infty \)
\end{itemize}

Given \( B, A \in \mathbb{R}^{n \times m} \) and an arbitrary \( s \in \mathbb{C} \), a matrix pencil is a family of matrices \( sB - A \), parametrized by \( s \). The pencil is called:

\begin{itemize}
    \item \textbf{Regular} when \( n = m \) and \( \det(sB - A) = p(s) \not\equiv 0 \).
    \item \textbf{Singular} when \( n \neq m \) or \( n = m \) and \( \det(sB - A) \equiv 0 \).
\end{itemize}

A regular pencil has \( p \) finite eigenvalues and \( q \) infinite eigenvalues with \( p + q = n \).

The equivalence theorem states:

\begin{itemize}
    \item A zero eigenvalue of \( sB - A \) is an infinite eigenvalue of \( B - sA \) and vice versa.
    \item A non-zero finite eigenvalue \( a_i \) of \( sB - A \) defines a non-zero finite eigenvalue \( \frac{1}{a_i} \) of \( B - sA \) and vice versa.
    \item An infinite eigenvalue of \( sB - A \) is a zero eigenvalue of \( B - sA \) and vice versa.
\end{itemize}

The dual GEP is critical for several reasons:

\begin{itemize}

        \item Helps in identifying potential points of instability in  systems
  
        \item Maintains sparsity and structure, reducing computational burden.
        \item Enables faster and more efficient simulations.

        \item Facilitates solving large-scale and complex eigenvalue problems.
        \item Useful for modern, interconnected power grids.
  
\end{itemize}

The dual GEP offers distinct advantages over the prime form, such as better sparsity preservation and complementary insights into system stability, making it a valuable tool for analyzing and optimizing complex systems.

\subsection{Applications and Extensions of the Continuous Electricity Market Model}

The continuous electricity market model provides a fundamental framework for understanding the dynamics of power markets. This section explores practical applications of the model, various methods for implementing it, and potential extensions to address real-world complexities.

\subsection*{Numerical Simulation}

Numerical simulation plays a crucial role in validating and exploring the dynamics predicted by the continuous electricity market model. Simulations help in understanding how the market responds to various conditions and disturbances.

\subsubsection*{Simulation Setup}

To perform numerical simulations, the differential equations governing the model are discretized using methods such as the Euler method or the Runge-Kutta method. The setup involves:

\begin{itemize}
    \item Discretizing time \( t \) into small intervals \( \Delta t \).
    \item Initializing the state variables \( \mathbf{S}(0) \), \( \mathbf{D}(0) \), and \( \lambda(0) \).
    \item Iteratively solving the differential equations for each time step.
\end{itemize}

\subsubsection*{Algorithm Implementation}

The simulation algorithm involves the following steps:

\begin{enumerate}
    \item \textbf{Initialization:} Set initial values for \( \mathbf{S}(0) \), \( \mathbf{D}(0) \), and \( \lambda(0) \).
    \item \textbf{Time Loop:} For each time step \( t_k \):
    \begin{itemize}
        \item Compute \( \mathbf{S}(t_k) \) and \( \mathbf{D}(t_k) \) using:
        \[
        \mathbf{S}(t_{k+1}) = \mathbf{S}(t_k) + \Delta t \cdot \boldsymbol{\alpha} (\lambda(t_k) \mathbf{1}_m - \mathbf{C})
        \]
        \[
        \mathbf{D}(t_{k+1}) = \mathbf{D}(t_k) + \Delta t \cdot \boldsymbol{\beta} (\mathbf{B} - \lambda(t_k) \mathbf{1}_n)
        \]
        \item Update the energy imbalance \( E(t_k) \):
        \[
        E(t_k) = \sum_{i=1}^{m} S_i(t_k) - \sum_{j=1}^{n} D_j(t_k)
        \]
        \item Adjust the power price \( \lambda(t_k) \):
        \[
        \lambda(t_{k+1}) = \lambda(t_k) + \Delta t \cdot (-k E(t_k) + h (\lambda_0 - \lambda(t_k)))
        \]
    \end{itemize}
    \item \textbf{Termination:} Stop when the simulation reaches the desired final time.
\end{enumerate}

\subsection*{Optimization Techniques}

Optimization techniques are employed to enhance the performance and efficiency of the electricity market. These methods can be used to optimize generator dispatch, minimize costs, and improve market stability.

\subsubsection*{Objective Functions}

Common objectives in the optimization of electricity markets include:

\begin{itemize}
    \item Minimizing the total generation cost:
    \[
    \text{Minimize} \quad \sum_{i=1}^{m} C_i S_i
    \]
    \item Maximizing social welfare:
    \[
    \text{Maximize} \quad \sum_{j=1}^{n} B_j D_j - \sum_{i=1}^{m} C_i S_i
    \]
\end{itemize}

\subsubsection*{Constraints}

Optimization problems are subject to various constraints such as:

\begin{itemize}
    \item Power balance:
    \[
    \sum_{i=1}^{m} S_i = \sum_{j=1}^{n} D_j
    \]
    \item Generation limits:
    \[
    S_i^{\text{min}} \leq S_i \leq S_i^{\text{max}}, \quad i = 1, \ldots, m
    \]
    \item Demand limits:
    \[
    D_j^{\text{min}} \leq D_j \leq D_j^{\text{max}}, \quad j = 1, \ldots, n
    \]
\end{itemize}

\subsubsection*{Solution Methods}

Various algorithms can be used to solve the optimization problems, including:

\begin{itemize}
    \item Linear Programming (LP)
    \item Mixed-Integer Linear Programming (MILP)
    \item Nonlinear Programming (NLP)
    \item Genetic Algorithms (GA)
    \item Particle Swarm Optimization (PSO)
\end{itemize}

\subsection*{Real-World Applications}

The continuous electricity market model can be extended and applied to address various real-world challenges in power systems.

\subsubsection*{Integration of Renewable Energy Sources}

Renewable energy sources such as wind and solar power introduce variability and uncertainty in the power supply. The model can be extended to include stochastic elements to represent these uncertainties and optimize the integration of renewables.

\subsubsection*{Demand Response Programs}

Demand response programs incentivize consumers to adjust their power usage in response to price signals. The model can incorporate demand response mechanisms to study their impact on market dynamics and stability.

\subsubsection*{Grid Stability and Resilience}

The model can be used to analyze the stability and resilience of the power grid under different scenarios, including faults, cyber-attacks, and natural disasters. This helps in designing strategies to enhance grid security and reliability.

\subsection*{Future Directions}

Future research can explore several directions to further enhance the continuous electricity market model:

\begin{itemize}
    \item Incorporating advanced machine learning algorithms for predictive analytics and optimization.
    \item Developing hybrid models that combine continuous and discrete elements to better capture the complexity of modern power systems.
    \item Enhancing the scalability of the model to handle large-scale power grids with high penetration of renewable energy sources.
    \item Integrating real-time data and control systems to create adaptive and responsive market mechanisms.
\end{itemize}

\section{Discrete Electricity Market Model}

Understanding the discrete dynamics of the electricity market is crucial for modeling the interactions between power producers and consumers. This section describes the discrete model that captures these dynamics through difference equations, focusing on the behavior of producers and consumers and the resulting market dynamics driven by energy imbalance.

\subsection*{Producer and Consumer Behavior}

In a market-driven environment, the behavior of power producers and consumers can be approximated by the following first-order difference equations, assuming network losses are negligible:

\begin{equation}
    S_i(k+1) = S_i(k) + \Delta t \cdot \alpha_i \left( \lambda(k) - (a_i + b_i S_i(k)) \right), \quad i = 1, \ldots, m
\end{equation}
\begin{equation}
    D_j(k+1) = D_j(k) + \Delta t \cdot \beta_j \left( (c_j - d_j D_j(k)) - \lambda(k) \right), \quad j = 1, \ldots, n
\end{equation}

where:
\begin{itemize}
    \item $S_i(k)$ is the power supply of producer $i$ at time step $k$.
    \item $D_j(k)$ is the power demand of consumer $j$ at time step $k$.
    \item $\lambda(k)$ is the price of power at time step $k$.
    \item $a_i$ is the fixed marginal cost of producer $i$.
    \item $b_i$ is the variable marginal cost coefficient of producer $i$.
    \item $c_j$ is the fixed marginal benefit of consumer $j$.
    \item $d_j$ is the variable marginal benefit coefficient of consumer $j$.
    \item $\alpha_i$ is the response rate of producer $i$.
    \item $\beta_j$ is the response rate of consumer $j$.
    \item $\Delta t$ is the discrete time step.
\end{itemize}

These equations describe the qualitative behavior where a generator increases its production when the price exceeds its marginal production cost, and loads act to increase consumption when marginal benefit exceeds the price. Each aims to maximize its profit or benefit by matching its marginal cost/benefit to the market price at equilibrium. The first-order difference equations allow approximate representation of generator ramp rates and lags in response to power price changes.

\subsection*{Energy Imbalance Market Dynamics}

In a synchronous power system, energy imbalance cannot be sustained indefinitely and must be reduced or driven to zero. Traditionally, this is managed by automatic generation control (AGC). In a real-time market-driven model, market mechanisms might fulfill this role by adjusting the power price based on energy imbalance. The energy imbalance-driven market dynamics can be represented by the following equations:

\begin{equation}
    E(k+1) = E(k) + \Delta t \cdot \left( \sum_{i=1}^{m} S_i(k) - \sum_{j=1}^{n} D_j(k) \right)
\end{equation}
\begin{equation}
    \lambda(k+1) = \lambda(k) + \Delta t \cdot (-k E(k) + h (\lambda_0 - \lambda(k)))
\end{equation}

where:
\begin{itemize}
    \item $E(k)$ is the system energy imbalance at time step $k$.
    \item $k$ is the power price response rate constant.
    \item $h$ is the market stabilizer gain for energy imbalance.
    \item $\lambda_0$ is the reference price of power.
    \item $\Delta t$ is the discrete time step.
\end{itemize}

The first equation represents the energy imbalance as the cumulative difference between the total power supply and the total power demand. The second equation adjusts the power price based on the energy imbalance and includes a feedback mechanism to stabilize the market price. The parameter $k$ dictates the sensitivity of price to the energy imbalance, and $h$ determines the rate at which the price returns to its reference value.

\subsection*{System Representation in Matrix Form}

To represent the system in a compact matrix form, we define the following vectors and matrices:

\begin{itemize}
    \item $\mathbf{S}(k) = [S_1(k), S_2(k), \ldots, S_m(k)]^T$ is the vector of power supplies from producers.
    \item $\mathbf{D}(k) = [D_1(k), D_2(k), \ldots, D_n(k)]^T$ is the vector of power demands from consumers.
    \item $\mathbf{C}(k) = [a_1 + b_1 S_1(k), a_2 + b_2 S_2(k), \ldots, a_m + b_m S_m(k)]^T$ is the vector of marginal costs for producers.
    \item $\mathbf{B}(k) = [c_1 - d_1 D_1(k), c_2 - d_2 D_2(k), \ldots, c_n - d_n D_n(k)]^T$ is the vector of marginal benefits for consumers.
    \item $\boldsymbol{\alpha} = \text{diag}(\alpha_1, \alpha_2, \ldots, \alpha_m)$ is the diagonal matrix of producer response rates.
    \item $\boldsymbol{\beta} = \text{diag}(\beta_1, \beta_2, \ldots, \beta_n)$ is the diagonal matrix of consumer response rates.
    \item $\mathbf{1}_m$ and $\mathbf{1}_n$ are vectors of ones of length $m$ and $n$, respectively.
\end{itemize}

The system of difference equations can then be written as:

\begin{equation}
    \mathbf{S}(k+1) = \mathbf{S}(k) + \Delta t \cdot \boldsymbol{\alpha} (\lambda(k) \mathbf{1}_m - \mathbf{C}(k))
\end{equation}
\begin{equation}
    \mathbf{D}(k+1) = \mathbf{D}(k) + \Delta t \cdot \boldsymbol{\beta} (\mathbf{B}(k) - \lambda(k) \mathbf{1}_n)
\end{equation}
\begin{equation}
    E(k+1) = E(k) + \Delta t \cdot (\mathbf{1}_m^T \mathbf{S}(k) - \mathbf{1}_n^T \mathbf{D}(k))
\end{equation}
\begin{equation}
    \lambda(k+1) = \lambda(k) + \Delta t \cdot (-k E(k) + h (\lambda_0 - \lambda(k)))
\end{equation}

This section provides a detailed description of the discrete electricity market model, including the producer and consumer behavior, energy imbalance market dynamics, and system representation in matrix form. This approach captures the discrete dynamics of power markets and can be used for practical implementations and analyses.

\subsection{Discrete Electricity Market Model with Memory}

Understanding the discrete dynamics of the electricity market is crucial for modeling the interactions between power producers and consumers. This section describes the discrete model that captures these dynamics through difference equations with memory, focusing on the behavior of producers and consumers and the resulting market dynamics driven by energy imbalance.

\subsection*{Producer and Consumer Behavior with Memory}

In a market-driven environment, the behavior of power producers and consumers can be approximated by the following first-order difference equations with memory, assuming network losses are negligible:

\begin{equation}
    S_i(k+1) = S_i(k) + \sum_{l=0}^{p} \alpha_i^l \left( \lambda(k-l) - (a_i + b_i S_i(k-l)) \right), \quad i = 1, \ldots, m
\end{equation}
\begin{equation}
    D_j(k+1) = D_j(k) + \sum_{l=0}^{p} \beta_j^l \left( (c_j - d_j D_j(k-l)) - \lambda(k-l) \right), \quad j = 1, \ldots, n
\end{equation}

where:
\begin{itemize}
    \item $S_i(k)$ is the power supply of producer $i$ at time step $k$.
    \item $D_j(k)$ is the power demand of consumer $j$ at time step $k$.
    \item $\lambda(k)$ is the price of power at time step $k$.
    \item $a_i$ is the fixed marginal cost of producer $i$.
    \item $b_i$ is the variable marginal cost coefficient of producer $i$.
    \item $c_j$ is the fixed marginal benefit of consumer $j$.
    \item $d_j$ is the variable marginal benefit coefficient of consumer $j$.
    \item $\alpha_i^l$ are the memory weights for producer $i$.
    \item $\beta_j^l$ are the memory weights for consumer $j$.
    \item $p$ is the memory length.
\end{itemize}

These equations describe the qualitative behavior where a generator increases its production when the price exceeds its marginal production cost, and loads act to increase consumption when marginal benefit exceeds the price. Each aims to maximize its profit or benefit by matching its marginal cost/benefit to the market price at equilibrium. The inclusion of memory allows the model to account for past behaviors and decisions.

\subsection*{Energy Imbalance Market Dynamics with Memory}

In a synchronous power system, energy imbalance cannot be sustained indefinitely and must be reduced or driven to zero. Traditionally, this is managed by automatic generation control (AGC). In a real-time market-driven model, market mechanisms might fulfill this role by adjusting the power price based on energy imbalance. The energy imbalance-driven market dynamics can be represented by the following equations with memory:

\begin{equation}
    E(k+1) = E(k) + \sum_{l=0}^{p} \left( \sum_{i=1}^{m} S_i(k-l) - \sum_{j=1}^{n} D_j(k-l) \right)
\end{equation}
\begin{equation}
    \lambda(k+1) = \lambda(k) + \sum_{l=0}^{p} \left( -k E(k-l) + h (\lambda_0 - \lambda(k-l)) \right)
\end{equation}

where:
\begin{itemize}
    \item $E(k)$ is the system energy imbalance at time step $k$.
    \item $k$ is the power price response rate constant.
    \item $h$ is the market stabilizer gain for energy imbalance.
    \item $\lambda_0$ is the reference price of power.
    \item $p$ is the memory length.
\end{itemize}

The first equation represents the cumulative energy imbalance as the sum of differences between the total power supply and the total power demand over the past $p$ time steps. The second equation adjusts the power price based on the energy imbalance and includes a feedback mechanism to stabilize the market price, considering the past $p$ time steps.

\subsection*{System Representation in Matrix Form with Memory}

To represent the system in a compact matrix form with memory, we define the following vectors and matrices:

\begin{itemize}
    \item $\mathbf{S}(k) = [S_1(k), S_2(k), \ldots, S_m(k)]^T$ is the vector of power supplies from producers.
    \item $\mathbf{D}(k) = [D_1(k), D_2(k), \ldots, D_n(k)]^T$ is the vector of power demands from consumers.
    \item $\mathbf{C}(k) = [a_1 + b_1 S_1(k), a_2 + b_2 S_2(k), \ldots, a_m + b_m S_m(k)]^T$ is the vector of marginal costs for producers.
    \item $\mathbf{B}(k) = [c_1 - d_1 D_1(k), c_2 - d_2 D_2(k), \ldots, c_n - d_n D_n(k)]^T$ is the vector of marginal benefits for consumers.
    \item $\boldsymbol{\alpha}^l = \text{diag}(\alpha_1^l, \alpha_2^l, \ldots, \alpha_m^l)$ is the diagonal matrix of producer response rates for memory step $l$.
    \item $\boldsymbol{\beta}^l = \text{diag}(\beta_1^l, \beta_2^l, \ldots, \beta_n^l)$ is the diagonal matrix of consumer response rates for memory step $l$.
    \item $\mathbf{1}_m$ and $\mathbf{1}_n$ are vectors of ones of length $m$ and $n$, respectively.
\end{itemize}

The system of difference equations with memory can then be written as:

\begin{equation}
    \mathbf{S}(k+1) = \mathbf{S}(k) + \sum_{l=0}^{p} \boldsymbol{\alpha}^l (\lambda(k-l) \mathbf{1}_m - \mathbf{C}(k-l))
\end{equation}
\begin{equation}
    \mathbf{D}(k+1) = \mathbf{D}(k) + \sum_{l=0}^{p} \boldsymbol{\beta}^l (\mathbf{B}(k-l) - \lambda(k-l) \mathbf{1}_n)
\end{equation}
\begin{equation}
    E(k+1) = E(k) + \sum_{l=0}^{p} (\mathbf{1}_m^T \mathbf{S}(k-l) - \mathbf{1}_n^T \mathbf{D}(k-l))
\end{equation}
\begin{equation}
    \lambda(k+1) = \lambda(k) + \sum_{l=0}^{p} \left( -k E(k-l) + h (\lambda_0 - \lambda(k-l)) \right)
\end{equation}

This section provides a detailed description of the discrete electricity market model with memory, including the producer and consumer behavior, energy imbalance market dynamics, and system representation in matrix form. This approach captures the discrete dynamics of power markets with historical influence and can be used for practical implementations and analyses.

\subsection{Importance of Memory and Discrete Modeling in Electricity Markets}

The dynamics of electricity markets are inherently complex, involving interactions between various stakeholders, including producers, consumers, and market operators. To effectively model these interactions, it is essential to capture not only the current state of the system but also its historical behavior. This is where the concepts of memory and discrete modeling become critically important.

\subsection*{Importance of Memory}

Memory in the context of electricity market modeling refers to the inclusion of past states and decisions in predicting future behavior. This concept is crucial for several reasons:

\begin{enumerate}
    \item \textbf{Historical Dependencies:} Power producers and consumers do not base their decisions solely on the current market prices or conditions. Their actions are influenced by past experiences and historical data, which can provide insights into market trends and patterns.
    \item \textbf{Smoother Transitions:} Incorporating memory helps in capturing the inertia and lag effects in the system. For instance, changes in production levels or consumption rates do not occur instantaneously but over a period influenced by previous states. Memory allows for a more realistic representation of these gradual changes.
    \item \textbf{Improved Predictive Accuracy:} Models that account for memory can better anticipate future market dynamics by recognizing patterns and anomalies in historical data. This can lead to more accurate forecasts and better decision-making.
    \item \textbf{Stability and Robustness:} Including memory in the model can contribute to the stability and robustness of the market. By considering past states, the system can be designed to avoid abrupt changes that might lead to instability or inefficiencies.
\end{enumerate}

\subsection*{Importance of Discrete Modeling}

Discrete modeling, as opposed to continuous modeling, involves representing the system at distinct time intervals. This approach offers several advantages:

\begin{enumerate}
    \item \textbf{Practical Implementation:} Electricity markets operate in discrete time intervals (e.g., hourly, daily). Discrete models align well with this operational reality, making them more practical for real-world applications.
    \item \textbf{Computational Efficiency:} Discrete models are often easier to implement and solve using numerical methods and algorithms. This can lead to faster computations, which is essential for real-time market operations and simulations.
    \item \textbf{Flexibility in Analysis:} Discrete models allow for the incorporation of complex behaviors and interactions that may be difficult to represent in a continuous framework. This flexibility is crucial for studying various scenarios and conducting sensitivity analyses.
    \item \textbf{Integration with Digital Systems:} Modern electricity markets rely heavily on digital systems for data collection, processing, and decision-making. Discrete models can be easily integrated with these systems, facilitating automated and efficient market operations.
\end{enumerate}

\subsection*{Future Studies: Stability and Computational Techniques}

Given the importance of memory and discrete modeling, future studies can focus on several exciting areas to further enhance our understanding and management of electricity markets:

\begin{enumerate}
    \item \textbf{Stability Analysis:} Investigating the stability of the electricity market models with memory can provide insights into the conditions that lead to stable or unstable market behavior. This includes studying the effects of various parameters and external disturbances on the market's stability.
    \item \textbf{Advanced Computational Techniques:} Employing advanced computational techniques such as machine learning, optimization algorithms, and parallel computing can improve the efficiency and accuracy of market simulations and forecasts. These techniques can handle large datasets and complex models, providing deeper insights into market dynamics.
    \item \textbf{Scenario Analysis and Forecasting:} Using discrete models with memory, researchers can conduct scenario analysis to evaluate the impact of different market policies, regulatory changes, and technological advancements. This can help in developing robust strategies for market management and planning.
    \item \textbf{Real-Time Market Operations:} Enhancing real-time market operations by incorporating adaptive algorithms that can quickly respond to changes in market conditions. This includes developing real-time pricing mechanisms, demand response strategies, and automated control systems.
    \item \textbf{Integration of Renewable Energy:} Studying the integration of renewable energy sources into the electricity market. Memory models can help in understanding the variability and intermittency of renewable energy and devising strategies to mitigate their impact on market stability.
\end{enumerate}

In conclusion, incorporating memory and discrete modeling into electricity market models provides a more realistic and practical representation of market dynamics. These enhancements can lead to improved predictive accuracy, stability, and efficiency in market operations. Future research in this area promises to unlock new possibilities for optimizing and managing electricity markets in an increasingly complex and dynamic environment.

\section{Continuous Electricity Market Model with Memory}

\subsection{Electricity Market Model with Delays}

In a market-driven environment, the behavior of power producers and consumers can be approximated by the following fractional differential equations with memory effects, assuming network losses are negligible:

\begin{equation}
    \frac{dS_i(t)}{dt} = \sum_{l=0}^{p} \alpha_i^l \left( \lambda(t-l\Delta t) - (a_i + b_i S_i(t-l\Delta t)) \right), \quad i = 1, \ldots, m
\end{equation}
\begin{equation}
    \frac{dD_j(t)}{dt} = \sum_{l=0}^{p} \beta_j^l \left( (c_j - d_j D_j(t-l\Delta t)) - \lambda(t-l\Delta t) \right), \quad j = 1, \ldots, n
\end{equation}

where:
\begin{itemize}
    \item $S_i(t)$ is the power supply of producer $i$ at time $t$.
    \item $D_j(t)$ is the power demand of consumer $j$ at time $t$.
    \item $\lambda(t)$ is the price of power at time $t$.
    \item $a_i$ is the fixed marginal cost of producer $i$.
    \item $b_i$ is the variable marginal cost coefficient of producer $i$.
    \item $c_j$ is the fixed marginal benefit of consumer $j$.
    \item $d_j$ is the variable marginal benefit coefficient of consumer $j$.
    \item $\alpha_i^l$ are the memory weights for producer $i$.
    \item $\beta_j^l$ are the memory weights for consumer $j$.
    \item $p$ is the memory length.
    \item $\Delta t$ is the discrete time step for memory.
\end{itemize}

These equations describe the qualitative behavior where a generator increases its production when the price exceeds its marginal production cost, and loads act to increase consumption when marginal benefit exceeds the price. Each aims to maximize its profit or benefit by matching its marginal cost/benefit to the market price at equilibrium. The inclusion of memory allows the model to account for past behaviors and decisions.

In a synchronous power system, energy imbalance cannot be sustained indefinitely and must be reduced or driven to zero. Traditionally, this is managed by automatic generation control (AGC). In a real-time market-driven model, market mechanisms might fulfill this role by adjusting the power price based on energy imbalance. The energy imbalance-driven market dynamics can be represented by the following equations with memory effects:

\begin{equation}
    \frac{dE(t)}{dt} = \sum_{i=1}^{m} S_i(t) - \sum_{j=1}^{n} D_j(t)
\end{equation}
\begin{equation}
    \frac{d\lambda(t)}{dt} = \sum_{l=0}^{p} \left( -k E(t-l\Delta t) + h (\lambda_0 - \lambda(t-l\Delta t)) \right)
\end{equation}

where:
\begin{itemize}
    \item $E(t)$ is the system energy imbalance at time $t$.
    \item $k$ is the power price response rate constant.
    \item $h$ is the market stabilizer gain for energy imbalance.
    \item $\lambda_0$ is the reference price of power.
    \item $p$ is the memory length.
    \item $\Delta t$ is the discrete time step for memory.
\end{itemize}

The first equation represents the energy imbalance as the difference between the total power supply and the total power demand. The second equation adjusts the power price based on the energy imbalance and includes a feedback mechanism to stabilize the market price, considering the past $p$ time steps.

To represent the system in a compact matrix form with memory, we define the following vectors and matrices:

\begin{itemize}
    \item $\mathbf{S}(t) = [S_1(t), S_2(t), \ldots, S_m(t)]^T$ is the vector of power supplies from producers.
    \item $\mathbf{D}(t) = [D_1(t), D_2(t), \ldots, D_n(t)]^T$ is the vector of power demands from consumers.
    \item $\mathbf{C}(t) = [a_1 + b_1 S_1(t), a_2 + b_2 S_2(t), \ldots, a_m + b_m S_m(t)]^T$ is the vector of marginal costs for producers.
    \item $\mathbf{B}(t) = [c_1 - d_1 D_1(t), c_2 - d_2 D_2(t), \ldots, c_n - d_n D_n(t)]^T$ is the vector of marginal benefits for consumers.
    \item $\boldsymbol{\alpha}^l = \text{diag}(\alpha_1^l, \alpha_2^l, \ldots, \alpha_m^l)$ is the diagonal matrix of producer response rates for memory step $l$.
    \item $\boldsymbol{\beta}^l = \text{diag}(\beta_1^l, \beta_2^l, \ldots, \beta_n^l)$ is the diagonal matrix of consumer response rates for memory step $l$.
    \item $\mathbf{1}_m$ and $\mathbf{1}_n$ are vectors of ones of length $m$ and $n$, respectively.
\end{itemize}

The system of differential equations with memory can then be written as:

\begin{equation}
    \frac{d\mathbf{S}(t)}{dt} = \sum_{l=0}^{p} \boldsymbol{\alpha}^l (\lambda(t-l\Delta t) \mathbf{1}_m - \mathbf{C}(t-l\Delta t))
\end{equation}
\begin{equation}
    \frac{d\mathbf{D}(t)}{dt} = \sum_{l=0}^{p} \boldsymbol{\beta}^l (\mathbf{B}(t-l\Delta t) - \lambda(t-l\Delta t) \mathbf{1}_n)
\end{equation}
\begin{equation}
    \frac{dE(t)}{dt} = \mathbf{1}_m^T \mathbf{S}(t) - \mathbf{1}_n^T \mathbf{D}(t)
\end{equation}
\begin{equation}
    \frac{d\lambda(t)}{dt} = \sum_{l=0}^{p} \left( -k E(t-l\Delta t) + h (\lambda_0 - \lambda(t-l\Delta t)) \right)
\end{equation}

This section provides a detailed description of the continuous electricity market model with memory, including the producer and consumer behavior, energy imbalance market dynamics, and system representation in matrix form. This approach captures the continuous dynamics of power markets with historical influence and can be used for practical implementations and analyses.

\subsection{Fractional Electricity Market Model}

Fractional calculus is a generalization of ordinary differentiation and integration to non-integer (fractional) orders. The Caputo fractional derivative of order \( \alpha \) of a function \( f(t) \) is defined as:

\begin{equation}
    {}^{C}D^{\alpha}_{t} f(t) = \frac{1}{\Gamma(n-\alpha)} \int_{0}^{t} (t-\tau)^{n-\alpha-1} f^{(n)}(\tau) \, d\tau
\end{equation}

where \( n = \lceil \alpha \rceil \) and \( \Gamma(\cdot) \) is the Gamma function.

In a market-driven environment, the behavior of power producers and consumers can be approximated by the following fractional-order differential equations, assuming network losses are negligible:

\begin{equation}
    {}^{C}D^{\alpha_i}_{t} S_i(t) = \lambda(t) - (a_i + b_i S_i(t)), \quad i = 1, \ldots, m
\end{equation}
\begin{equation}
    {}^{C}D^{\beta_j}_{t} D_j(t) = (c_j - d_j D_j(t)) - \lambda(t), \quad j = 1, \ldots, n
\end{equation}

where:
\begin{itemize}
    \item $S_i(t)$ is the power supply of producer $i$ at time $t$.
    \item $D_j(t)$ is the power demand of consumer $j$ at time $t$.
    \item $\lambda(t)$ is the price of power at time $t$.
    \item $a_i$ is the fixed marginal cost of producer $i$.
    \item $b_i$ is the variable marginal cost coefficient of producer $i$.
    \item $c_j$ is the fixed marginal benefit of consumer $j$.
    \item $d_j$ is the variable marginal benefit coefficient of consumer $j$.
    \item ${}^{C}D^{\alpha_i}_{t}$ denotes the Caputo fractional derivative of order $\alpha_i$ for producer $i$.
    \item ${}^{C}D^{\beta_j}_{t}$ denotes the Caputo fractional derivative of order $\beta_j$ for consumer $j$.
\end{itemize}

These equations describe the behavior of producers and consumers with memory effects introduced through the fractional derivatives. This allows for a more accurate modeling of real-world dynamics, where past states influence current behavior.

In a synchronous power system, energy imbalance cannot be sustained indefinitely and must be reduced or driven to zero. Traditionally, this is managed by automatic generation control (AGC). In a real-time market-driven model, market mechanisms might fulfill this role by adjusting the power price based on energy imbalance. The energy imbalance-driven market dynamics can be represented by the following fractional-order equations:

\begin{equation}
    {}^{C}D^{\gamma}_{t} \lambda(t) = -H_d \lambda(t) + K E (\omega_{ref} - \omega_{CoI}(t))
\end{equation}
\begin{equation}
    \frac{dE(t)}{dt} = \sum_{i=1}^{m} S_i(t) - \sum_{j=1}^{n} D_j(t)
\end{equation}

where:
\begin{itemize}
    \item $H_d$ is the damping coefficient.
    \item $K E (\omega_{ref} - \omega_{CoI}(t))$ is the control input based on the frequency deviation.
    \item ${}^{C}D^{\gamma}_{t}$ denotes the Caputo fractional derivative of order $\gamma$.
    \item $E(t)$ is the system energy imbalance at time $t$.
\end{itemize}

The first equation represents the fractional-order dynamics of the power price adjustment based on energy imbalance. The second equation represents the energy imbalance as the difference between the total power supply and the total power demand.

To represent the system in a compact matrix form with fractional derivatives, we define the following vectors and matrices:

\begin{itemize}
    \item $\mathbf{S}(t) = [S_1(t), S_2(t), \ldots, S_m(t)]^T$ is the vector of power supplies from producers.
    \item $\mathbf{D}(t) = [D_1(t), D_2(t), \ldots, D_n(t)]^T$ is the vector of power demands from consumers.
    \item $\mathbf{C}_{gi} = [a_1 + b_1 S_1(t), a_2 + b_2 S_2(t), \ldots, a_m + b_m S_m(t)]^T$ is the vector of marginal costs for producers.
    \item $\mathbf{C}_{di} = [c_1 - d_1 D_1(t), c_2 - d_2 D_2(t), \ldots, c_n - d_n D_n(t)]^T$ is the vector of marginal benefits for consumers.
    \item $\boldsymbol{\alpha} = \text{diag}(\alpha_1, \alpha_2, \ldots, \alpha_m)$ is the diagonal matrix of producer response rates.
    \item $\boldsymbol{\beta} = \text{diag}(\beta_1, \beta_2, \ldots, \beta_n)$ is the diagonal matrix of consumer response rates.
    \item $\mathbf{1}_m$ and $\mathbf{1}_n$ are vectors of ones of length $m$ and $n$, respectively.
\end{itemize}

The system of fractional-order differential equations can then be written as:

\begin{equation}
    \boldsymbol{\alpha} \, {}^{C}D^{\alpha}_{t} \mathbf{S}(t) = \lambda(t) \mathbf{1}_m - \mathbf{C}_{gi}
\end{equation}
\begin{equation}
    \boldsymbol{\beta} \, {}^{C}D^{\beta}_{t} \mathbf{D}(t) = -\lambda(t) \mathbf{1}_n + \mathbf{C}_{di}
\end{equation}
\begin{equation}
    {}^{C}D^{\gamma}_{t} \lambda(t) = -H_d \lambda(t) + K E (\omega_{ref} - \omega_{CoI}(t))
\end{equation}
\begin{equation}
    \frac{dE(t)}{dt} = \mathbf{1}_m^T \mathbf{S}(t) - \mathbf{1}_n^T \mathbf{D}(t)
\end{equation}

\subsection*{Importance of Fractional-Order Modeling and Machine Learning}

Incorporating fractional-order derivatives into the electricity market model introduces memory effects, reflecting the impact of past states on current dynamics. This approach enhances the model's realism and predictive capability, capturing the conservative behavior of market participants as they respond to historical data.

Additionally, the use of fractional-order models can be complemented with machine learning techniques to further enhance their predictive power and efficiency. Machine learning algorithms can be employed to:

\begin{itemize}
    \item \textbf{Optimize Model Parameters:} Machine learning can be used to optimize the parameters of the fractional-order model, ensuring that it accurately reflects observed market behavior.
    \item \textbf{Predict Future States:} By training on historical data, machine learning models can predict future market states, providing valuable insights for decision-making.
    \item \textbf{Identify Patterns:} Advanced algorithms can identify complex patterns and correlations in market data that might not be apparent through traditional analysis.
    \item \textbf{Enhance Stability Analysis:} Machine learning techniques can be used to analyze the stability of the fractional-order model, identifying conditions that lead to stable or unstable market behavior.
\end{itemize}

\section{Conclusion and Future Work}

In this paper, we presented the development and analysis of three distinct electricity market models: the continuous model, the discrete model, and the fractional-order model. Each model addresses unique challenges and characteristics of the electricity market, providing valuable insights into its dynamics and behavior.

\subsection*{Summary of Models}

\textbf{Continuous Model:} This model employs continuous mathematical techniques to represent the smooth and ongoing changes in electricity generation, transmission, distribution, and consumption. It is particularly useful for understanding long-term trends and steady-state behaviors. The continuous model captures the interactions between power producers and consumers through differential equations and represents energy imbalance dynamics with feedback mechanisms.

\textbf{Discrete Model:} The discrete model simulates the electricity market using discrete time intervals, making it ideal for analyzing scenarios where market events occur at specific points in time, such as bidding processes and market clearing. This model provides insights into short-term fluctuations and the discrete nature of market operations. By incorporating memory effects, the discrete model can account for historical behaviors and decisions, offering a more realistic representation of market dynamics.

\textbf{Fractional-Order Model:} Introducing fractional calculus into the electricity market model allows for capturing memory effects and hereditary properties of the market. Fractional-order models provide a nuanced understanding of complex temporal dependencies and long-term correlations that are not easily addressed by traditional models. This approach enhances the model's realism and predictive capability, reflecting the conservative behavior of market participants as they respond to historical data.

\subsection*{Importance of Memory and Discrete Modeling}

Incorporating memory and discrete modeling into electricity market models is crucial for capturing the historical dependencies and lag effects in the system. Memory allows for a more accurate modeling of real-world dynamics, where past states influence current behavior. Discrete modeling aligns well with the operational reality of electricity markets, which operate in discrete time intervals, and facilitates practical implementation and computational efficiency.

\subsection*{Future Work}

Future research can explore several exciting areas to further enhance the understanding and management of electricity markets:

\textbf{Stability Analysis:} Investigating the stability of the electricity market models with memory can provide insights into the conditions that lead to stable or unstable market behavior. This includes studying the effects of various parameters and external disturbances on the market's stability.

\textbf{Advanced Computational Techniques:} Employing advanced computational techniques such as machine learning, optimization algorithms, and parallel computing can improve the efficiency and accuracy of market simulations and forecasts. These techniques can handle large datasets and complex models, providing deeper insights into market dynamics.

\textbf{Scenario Analysis and Forecasting:} Using discrete models with memory, researchers can conduct scenario analysis to evaluate the impact of different market policies, regulatory changes, and technological advancements. This can help in developing robust strategies for market management and planning.

\textbf{Real-Time Market Operations:} Enhancing real-time market operations by incorporating adaptive algorithms that can quickly respond to changes in market conditions. This includes developing real-time pricing mechanisms, demand response strategies, and automated control systems.

\textbf{Integration of Renewable Energy:} Studying the integration of renewable energy sources into the electricity market. Memory models can help in understanding the variability and intermittency of renewable energy and devising strategies to mitigate their impact on market stability.

In conclusion, the comprehensive approach of combining continuous, discrete, and fractional-order models provides a robust framework for understanding and analyzing the dynamics of electricity markets. These models offer valuable insights into market behavior, stability, and efficiency, and future research will further enhance their applicability and effectiveness in the evolving energy landscape.



\section*{Competing of Interest}

The author declare no conflict of interest.



\end{document}